\documentclass{article}
\usepackage[english]{babel}
\usepackage[utf8]{inputenc}
\usepackage{johd}

\usepackage{url}
\usepackage{tikz,pgfplots}
\usetikzlibrary{positioning}
\pgfplotsset{compat=1.14} 
\tikzstyle{vertex} = [fill,shape=circle,node distance=30pt]
\tikzstyle{edge} = [fill,opacity=.6,fill opacity=.5,line cap=round, line join=round, line width=10pt]
\tikzstyle{elabel} =  [fill,shape=circle,node distance=30pt,fill opacity=.9]
\definecolor{mygray}{gray}{0.95}
\pgfdeclarelayer{background}
\pgfsetlayers{background,main}
\usepackage[most]{tcolorbox}
\definecolor{mypurple}{rgb}{0.59, 0.44, 0.84}
\usepackage{hyperref}

\usepackage[utf8]{inputenc}
\usepackage{amsmath,amssymb}%amsthm
\usepackage{color}
\usepackage{mathrsfs}
\usepackage{eucal}
\usepackage{mathtools}
\usepackage{nameref}        % reference section by name
\usepackage{algorithmic}
\usepackage{algorithm}
\usepackage{pgfplots}
\pgfplotsset{compat=1.16}

\newtheorem{thm}{Theorem}
\newtheorem{remark}{Remark}

\newtheorem{prop}{Property}

\newtheorem{defn}{Definition}

\newcommand{\R}{\mathbb{R}}

\newcommand{\Z}{\mathbb{Z}}

\newcommand{\g}{\mathcal{G}}

\newcommand{\Vs}{\mathcal{V}}
\newcommand{\Es}{\mathcal{E}}
\newcommand{\xv}{\mathbf{x}}
\newcommand{\uv}{\mathbf{u}}

\newcommand{\Av}{\mathbf{A}}
\newcommand{\Cv}{\mathbf{C}}
\newcommand{\Bv}{\mathbf{B}}

\newcommand{\tA}{\textsf{A}}
\newcommand{\tT}{\textsf{T}}

\newcommand{\yv}{\mathbf{y}}
\newcommand{\hv}{\mathbf{h}}
\newcommand{\fv}{\mathbf{f}}
\newcommand{\gv}{\mathbf{g}}
\newcommand{\Id}{\mathbf{I}}
\newcommand{\tra}{\prime}
\newcommand{\h}{\mathcal{H}}

\title{\bf Geometric Observability of Hypergraphs}

\author{Joshua Pickard$^{a}$$^{*}$, Cooper Stansbury$^{a}$, Amit Surana$^{b}$,\\ Indika Rajapakse$^{ac}$, and Anthony Bloch$^{c}$ \\
        \small $^{a}$Department of Computational Medicine and Bioinformatics, University of Michigan, Ann Arbor MI, USA \\
        \small $^{b}$RTX Technology Research Center, East Hartford, CT USA \\
        \small $^{c}$Department of Mathematics, University of Michigan, Ann Arbor MI, USA \\\\
}
\date{\today}

\begin{document}
\maketitle
\begin{abstract}
In this paper we consider aspects of geometric observability for hypergraphs, extending our earlier work from the uniform to the nonuniform case. 
Hypergraphs, a generalization of graphs, allow hyperedges to connect multiple nodes and unambiguously represent multi-way relationships which are ubiquitous in many real-world networks
including those that arise in biology. We consider polynomial dynamical systems with linear outputs defined according to hypergraph structure, and we propose methods to evaluate local, weak observability.

\vspace{5mm}

\end{abstract}

%%%%%%%%%%%%%%%%%%%%%%%%%%%%%%%%%%%%%%%%%%%%%%%%%%%%%%%%%%%%%%%%%%%%%%%%%%%%%%%%
\section{INTRODUCTION}
% In this paper, we develop a framework for studying observability of complex networks represented as hypergraphs.
% The ability to monitor, predict, and control complex, networked systems is a fundamental and crucial task with widespread applications in various domains, including social/communications systems, life sciences, and security/defense, among others others.
We consider here geometric aspects of observability with applications to hypergraphs % and % \textcolor{blue}{unless we put mouse example, we should probably not emphasize biological. Tony: added in "motivated by"} 
motivated by  biological networks. 
%\textcolor{blue}{Hypergraphs, a generalization of graphs, allow hyperedges to connect multiple nodes and unambiguously represent multi-way relationships which are ubiquitous in many real-world networks \cite{battiston2020networks}.
In particular, we consider
nonuniform hypergraphs in which hyperedges can contain arbitrary finite number of nodes.  This extends our work in \cite{PiSuBlRa2023} which was restricted to observability analysis of uniform hypergraphs.

Hypergraphs extend classical graph theory by allowing hyperedges to have more than two vertices. 
This is important because may complex systems cannot be described by pairwise interactions only. 
Examples include social networks as well as complex physical and biological networks such as polymers
and proteins. In \cite{chen2021controllability}  we consider a mouse neuron model
as well as interaction in the human genome which can be captured by hypergraph models. Allowing hypergraphs to be nonuniform allows one to capture the dynamics of 
complex asymmetrical systems. 

Two fundamental questions arise when considering nonuniform hypergraph observability:
\begin{itemize}
    \item (Q1) Is a set of sensor nodes sufficient to render a network observable?
    \item (Q2) How do higher order interactions contribute to observability?
\end{itemize}

Observability of networks has been considered from an array of perspectives; see \cite{montanari2020observability} and references therein. To address Q1, structural observability utilizes network topology \cite{lin1974structural}; dynamic observability applies classical matrix properties, particle filtering \cite{montanari2019particle}, or the observability Gramian \cite{summers2014optimal}; and topological observability explores the relationship between observability and graph topologies, \cite{liu2013observability,su2017analysis}. Despite the array of approaches, hypergraph observability and Q2 remain relatively unexplored.

To investigate Q2, this paper makes the following contributions:
\begin{itemize}
    \item We propose a multitensor representation of the adjacency structure and dynamics of nonuniform hypergraphs.
    \item We develop a nonlinear observability test for nonuniform hypergraphs and demonstrate it on several canonical hypergraph topologies.
\end{itemize}
We focus on the concept of weak local observability for polynomial systems, and to overcome the limitations of local observability, our computations leverage symbolic calculations to offer a global observability test.

We begin with necessary background on nonlinear observability (Section 2) and uniform hypergraphs (Section 3). We extend these results to non-uniform hypergraphs, where we develop our calculations for the observability of nonuniform hypergraphs (Section 4).

\section{NONLINEAR OBSERVABILITY BACKGROUND}\label{sec:obs}
Nonlinear controllability and observability were introduced in the seminal work \cite{hermankrener};
see the work of \cite{baillieul1981controllability} for polynomial dynamics. In contrast to linear systems, nonlinear observability, which is based similarly on the distinguishability of system states, has multiple varieties, such as local, weak, and global observability, \cite{hermankrener,sontag1984concept,gerbet2020global}. Unfortunately, unlike the linear case where the Kalman rank condition or Popov-Belevitch-Hautus test simply  determine observability, no easy criteria exist for nonlinear systems.

Consider the affine control system $\Sigma$,
\begin{equation*}
\Sigma\begin{cases}
\dot{\xv}=\fv(\xv,\uv)=\hv_0(\xv)+\sum_{i=1}^k\hv_i(\xv)u_i\\
\yv=\gv(\xv)
\end{cases}
\end{equation*}
where, $\uv=(u_1,\dots,u_k)^\tra\in \R^k$ denotes the input vector, $\xv\in M\subset \R^n$ is the state vector and $\yv\in\R^m$ is the output/measurement vector. We assume that $\Sigma$ is analytic, i.e., the functions $\hv_i:M\rightarrow M, i=0,\dots,k$ and $g_i:M\rightarrow \R, i=1,\dots,m$ where $\gv=(g_1,\dots,g_m)^\tra$ are assumed to be analytic functions defined on $M$. We also have to assume $\Sigma$ is complete, that is, for every bounded measurable input $\uv(t)$ and every $\xv_0\in M$ there exists a solution $\xv(t)$ of $\Sigma$ such that $\xv(0) = \xv_0$ and $\xv(t)\in M$ for all $t\in \R$. We review different notions of observability from \cite{obsThesis} which are equivalent to those introduced in \cite{hermankrener}, but use a slightly different terminology.

\begin{defn}
Let $U$ be an open subset of $M$. A pair of points $\xv_0$ and $\xv_1$ in $M$ are called $U$-\textit{distinguishable } if there exists a measurable input $\uv(t)$ defined on the interval $[0,T]$ that generates solutions $\xv_0(t)$ and $\xv_1(t)$ of system $\Sigma$ satisfying $\xv_i(0)=\xv_i,i=0,1$ such that $\xv_i(t)\in U$ for $t\in[0,T]$ and $\gv(\xv_0(t))\neq \gv(\xv_1(t))$ for some $t\in [0,T]$. We denote by $I(\xv_0,U)$ all points $\xv_1\in U$ that are not $U$--\textit{distinguishable } from $\xv_0$
\end{defn}

\begin{defn}
The system $\Sigma$ is \textit{observable} at $\xv_0\in M$ if $I(\xv_0,M)=\xv_0$.
\end{defn}

\begin{defn}
The system $\Sigma$ is \textit{locally observable} at $\xv_0\in M$ if for every open neighbourhood $U$ of $\xv_0$, $I(\xv_0,U)=\xv_0$
\end{defn}

Local observability implies observability. On the other hand, since $U$ can be chosen arbitrarily small, local observability implies that we can distinguish between neighboring points instantaneously. Both the definitions above ensure that a point $\xv_0\in M$   can be distinguished from every other point in $M$. It is often sufficient to distinguish between neighbours in $M$, which leads to the following two notions of observability.

\begin{defn}
The system $\Sigma$ is \textit{weakly observable} at $\xv_0\in M$ if $\xv_0$ has an open neighbourhood $U$ such that $I(\xv_0,M)\bigcap U=\xv_0$.
\end{defn}

\begin{defn}
The system $\Sigma$ is \textit{locally weakly observable} at $\xv_0\in M$ if $\xv_0$ has an open neighbourhood $U$ such that for every open neighbourhood $V$ of $\xv_0$ contained in $U$, $I(\xv_0,V)=\xv_0$.
\end{defn}

As we can set $U=M$, local observability implies local weak observability. The local weakly observability lends itself to a simple algebraic test.  Let $\h$ be the \textit{observation space},
\begin{eqnarray*}
\h&=&\{L_{\hv_{i_1}}L_{\hv_{i_2}}\cdots L_{\hv_{i_r}}(g_{i}): r\geq 0, i_j=0,\cdots,k,\notag\\
 &&i=1,\cdots,m\} \label{eq:defG}
\end{eqnarray*}
{where $L_f$ denotes the Lie derivative w.r.t. vector field $f$,} and let 
\begin{equation*}\label{eq:defdG}
d\h=\mbox{span}_{\R_x}\{d\phi:\phi \in \h\}
\end{equation*}
be the space spanned by the gradients of the elements of $\h$, where $\R_x$ is the space of meromorphic functions on $M$. The following result was proved in \cite{hermankrener}, see Theorems $3.1$ and $3.11$.
\begin{thm}\label{thm1}
The analytic system $\Sigma$ is locally weakly observable for all $\xv$ in an open dense set of $M$ if and only if $dim_{\R_x}(d\h) = n$.
\end{thm}
\begin{remark}
Here $\mbox{dim}_{\R_x}(d\h)$ is the generic or maximal rank of $d\h$, that is, $\mbox{dim}_{\R_x}(d\h)= \max_{\xv\in M}(dim_{\R} d\h(\xv))$.
\end{remark}
For system $\Sigma$ with no control inputs, i.e. $\hv_i\equiv0,i=1,\cdots,k$ the condition for local weak observability simplifies to checking,
\begin{equation}\label{eq:obsr}
\mbox{rank}(\mathcal{O}(\xv))=n,
\end{equation}
where $\mathcal{O}(\xv)$ is the nonlinear observability matrix (NOM),
\begin{equation}\label{eq:obsrmat}
    \mathcal{O}(\xv)=\nabla_{\xv}\begin{pmatrix}
    L^0_{\hv_0}\gv(\xv)\\
    L^1_{\hv_0}\gv(\xv)\\
    \vdots\\
    L^r_{\hv_0}\gv(\xv)
    \end{pmatrix},
\end{equation}
for some $r\in \Z$. One can use symbolic computation to check the generic rank condition (\ref{eq:obsr}) as performed by Sedoglavic's algorithm, \cite{sedoglavic2001probabilistic}.

\begin{remark}\label{remrn}
In general the value of $r$  to use in (\ref{eq:obsrmat}) is not known apriori. For an analytic system $\Sigma$, $r$ can be set to the state dimension $n$, see Theorem $4.1$ in \cite{obsThesis}.
\end{remark}

\begin{remark}\label{algeom}
For a polynomial system $\Sigma$, observability has also been studied from the perspective of algebraic geometry, see \cite{gerbet2020global} and references therein.
\end{remark}

% \begin{remark}
We adopt the use of local weak observability as the notion of nonlinear observability throughout the remainder of this paper.
% \end{remark}

\section{UNIFORM HYPERGRAPHS}\label{sec: uniform hypergraphs}
A \textit{undirected hypergraph} is given by  $\g=\{\Vs,\Es\}$ where $\Vs$ is a finite set and $\Es\subseteq \mathcal{P}(\Vs)\setminus \{\emptyset\}$, the power set of $\Vs$ (i.e. the set of all subsets of $\Vs$). The elements of $\Vs$ are called the nodes, and the elements of $\Es$ are called the hyperedges. Cardinality $|e|$ of a hyperedge $e\in \Es$  is number of nodes contained in it. We will consider hypergraphs with no self-loops, so that all hyperedges have $|e|\geq 2$.   A hypergraph is $k$-uniform if all hyperedges contain exactly $k$ vertices.
\subsection{Uniform Hypergraph Structure}
% We note that in this definition of hypergraph we do not allow for repeated nodes within an hyperedge (often called hyperloops).  The degree $d(v)$ of a node $v\in\Vs$ is $d(v)=|\{e\in \Es|v\in e\}|$, where $|\cdot|$ denotes set cardinality. The degree of an hyperedge $e$ is denoted by $d(e) = |e|$. For $k$-uniform hypergraphs, the degree of each hyperedge is the same, i.e. $d(e) = k$.

\begin{defn}
Let $\g= \{\Vs,\Es\}$ be a $k$-uniform hypergraph with $n=|\Vs|$ nodes. The \textit{adjacency tensor} $\tA\in\mathbb{R}^{n\times n\times \dots\times n}$ of $\g$ is a $k$-th order, $n$-dimensional, supersymmetric tensor, which is defined as
\begin{equation}\label{eq:98}
    \tA_{j_1j_2\dots j_k} = 
    \begin{cases}
        \frac{1}{(k-1)!} & \text{if $\{j_1,j_2,\dots,j_k\}\in \Es$}\\
        0 &  \text{otherwise}\end{cases}.
\end{equation}

% Similarly to standard graphs, the \textit{degree} of vertex $j$ of a uniform hypergraph is defined as
% \begin{equation}\label{eq:r6}
% d_j=\sum_{j_2=1}^n\sum_{j_3=1}^n\dots \sum_{j_k=1}^n\tA_{jj_2j_3\dots j_k}.
% \end{equation}
\end{defn}

% Note that the choice of the nonzero coefficient $\frac{1}{(k-1)!}$ in (\ref{eq:98}) guarantees that the degree of each node is equal to the number of hyperedges that contain that node, which is consistent with the notion of degree in standard graphs. The \textit{degree distribution} of a hypergraph is the distribution of the degrees over all nodes. If all nodes have the same degree $d$, then $\g$ is called \textit{$d$-regular}.

We recall definitions of uniform hypergraph chain, ring, star and complete hypergraphs following \cite{chen2021controllability}.

\begin{defn}\label{def: hyperchain}
 A $k$-uniform \textit{hyperchain} is a sequence of $n$ nodes such that every $k$ consecutive nodes are adjacent, i.e., nodes $j,j+1,\dots,j+k-1$ are contained in one hyperedge for $j=1,2,\dots,n-k+1$.
\end{defn}

\begin{defn}
A $k$-uniform \textit{hyperring} is a sequence of $n$ nodes such that every $k$ consecutive nodes are adjacent, i.e., nodes $\sigma_n(j),\sigma_n(j+1),\dots,\sigma_n(j+k-1)$ are contained in one hyperedge for $j=1,2,\dots,n$, where $\sigma_n(j)=j$ for $j\leq n$ and $\sigma_n(j)=j-n$ for $j>n$.
\end{defn}

\begin{defn}
A $k$-uniform \textit{hyperstar} is a collection of $k-1$ internal nodes that are contained in all the hyperedges, and $n-k+1$ leaf nodes such that every leaf node is contained in one hyperedge with the internal nodes.
\end{defn}

\begin{defn}\label{def: complete hypergraph}
A $k-$uniform \textit{complete hypergraph} is a set of $n$ vertices with all $\binom{n}{k}$ possible hyperedges.
\end{defn}

Note that when $k=2$, definitions \ref{def: hyperchain} - \ref{def: complete hypergraph} are reduced to standard chains, rings, stars and complete graphs.

\subsection{Uniform Hypergraph Dynamics with Linear Outputs}
We consider the homogeneous polynomial/multilinear time-invariant dynamics on a $k$-uniform hypergraph with linear system outputs as introduced in \cite{PiSuBlRa2023}.

% Following \cite{chen2021controllability} we represent the dynamics of a uniform hypergraph by homogeneous polynomial/multilinear time-invariant differential equations with linear outputs as follows.

\begin{defn}
Given a $k$-uniform undirected hypergraph $\g$ with $n$ nodes, the dynamics of $\g$ with outputs $\yv\in \R^m$ is defined as
\begin{equation}\label{eq:hypdyn1}
    \g\begin{cases}
        \dot{\xv} & = \fv(\xv)=\tA\xv^{k-1}\\
        \yv & = \gv(\xv)=\Cv\xv,        
    \end{cases}
\end{equation}
where $\tA\in\mathbb{R}^{n\times n\times \dots \times n}$ is the adjacency tensor of $\g$, $\Cv\in\mathbb{R}^{m\times n}$ is the output matrix, and $\tA\xv^{k-1}$ is
% \textcolor{blue}{
\begin{equation*}
\tA\xv^{k-1}=\tA\times_2\xv\times_3\xv\cdots\times_k\xv,
\end{equation*}
% }
with $\times$ being the Tucker product, see \cite{PiSuBlRa2023} for details.
\end{defn}

 The \textit{matrix tensor multiplication} $\textsf{X} \times_{n} \textbf{A}$ along mode $n$ for a matrix $\textbf{A}\in  \mathbb{R}^{I\times J_n}$ is defined by
$
(\textsf{X} \times_{n} \textbf{A})_{j_1j_2\dots j_{n-1}ij_{n+1}\dots j_N}=\sum_{j_n=1}^{J_n}\textsf{X}_{j_1j_2\dots j_n\dots j_N}\textbf{A}_{ij_n}.
$
This product can be generalized to what is known as the \textit{Tucker product}, for $\textbf{A}_n\in \mathbb{R}^{I_n\times J_n}$,
\begin{align*}\label{eq5}
&\textsf{X}\times_1 \textbf{A}_1 \times_2\textbf{A}_2\times_3\dots \times_{N}\textbf{A}_N\\
&=\textsf{X}\times\{\textbf{A}_1,\textbf{A}_2,\dots,\textbf{A}_N\}\in  \mathbb{R}^{I_1\times I_2\times\dots \times I_N}.
\end{align*}

\begin{figure}[t!]
\centering
\tcbox[colback=white,top=5pt,left=5pt,right=-5pt,bottom=5pt]{
\begin{tikzpicture}[scale=0.85, transform shape]
\node[vertex,text=white,scale=0.7] (v1) {1};
\node[vertex,above right = 11pt and 1.5pt of v1,text=white,scale=0.7] (v2) {2};
\node[vertex,below right = 11pt and 1.5pt of v2,text=white,scale=0.7] (v3) {3};

\node[vertex,below right = 3pt and 40pt of v2,white,scale=0.5] (v7){};
\node[vertex,below right = 3pt and 6pt of v2,white,scale=0.5] (v8){};
\node[vertex,below right = -9pt and 12pt of v2,white,scale=0.7,text=black] (v9){  $\dot{\textbf{x}}$=$\textbf{A}\textbf{x}$};
\node[rectangle,below right = -20pt and 40pt of v2,white,scale=0.7,text=black] (v6) {
$\Large
\begin{cases}
\dot{x}_1 = x_2+x_3\\
\dot{x}_2 = x_1 + x_3\\
\dot{x}_3 = x_1 + x_2
\end{cases}
$};
\node[vertex,above of=v2,node distance=20pt,white,text=black,scale=0.7] (text1) {\Large\textbf{A}};
\path [->,shorten >=1pt,shorten <=1pt, thick](v8) edge node[left] {} (v7);

\node[vertex,right of=v3,node distance=120pt,text=white,scale=0.7](w1) {1};
\node[vertex,above right = 11pt and 1.5pt of w1,text=white,scale=0.7] (w2) {2};
\node[vertex,below right = 11pt and 1.5pt of w2,text=white,scale=0.7] (w3) {3};
\node[vertex,below right = 3pt and 40pt of w2,white,scale=0.5] (w7){};
\node[vertex,below right = 3pt and 6pt of w2,white,scale=0.5] (w8){};
\node[vertex,below right = -10pt and 12pt of w2,white,scale=0.7,text=black] (w9){  $\dot{\textbf{x}}$=$\textsf{A}\textbf{x}^{2}$};

\node[rectangle,below right = -20pt and 40pt of w2,white,scale=0.7,text=black] (w6) {
$\Large
\begin{cases}
\dot{x}_1 = x_2x_3\\
\dot{x}_2 = x_1x_3\\
\dot{x}_3 = x_1x_2
\end{cases}
$};
\node[vertex,above of=w2,node distance=20pt,white,text=black,scale=0.7] (text1) {\Large\textbf{B}};

\path [->,shorten >=1pt,shorten <=1pt, thick](w8) edge node[left] {} (w7);

\begin{pgfonlayer}{background}
\draw[edge,color=orange,line width=8pt] (v1) -- (v2);
\draw[edge,color=red,line width=8pt] (v2) -- (v3);
\draw[edge,color=green,line width=8pt] (v3) -- (v1);

\draw[edge,color=orange,line width=8pt] (w1) -- (w2) -- (w3) -- (w1);
\end{pgfonlayer}

\end{tikzpicture}
}
\caption{Graphs versus uniform hypergraphs. (A) Standard graph with three nodes and edges $e_1=\{1,2\}$, $e_2=\{2,3\}$ and $e_3=\{1,3\}$, and its corresponding linear dynamics. (B) 3-uniform hypergraph with three nodes and a hyperedge $e_1=\{1,2,3\}$, and its corresponding nonlinear dynamics.}
\label{fig:1}
\end{figure}
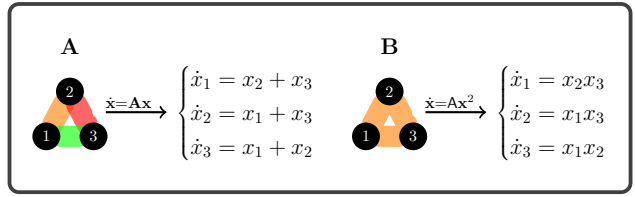

See Fig. \ref{fig:1} for an example of uniform hypergraph and associated dynamics. All the interactions are characterized using multiplications instead of the additions that are typically used in a standard graph based representation. For detailed discussion on relationship between graph vs. hypergraph dynamic representation, see \cite{chen2021controllability, pickard2023kronecker}.

\section{NON-UNIFORM HYPERGRAPHS}\label{sec: obsv uniform hypergraphs}

In this section, we recast the homogeneous polynomial/multilinear hypergraph dynamical system in terms of the Kronecker product and derive a construction of the corresponding nonlinear observability matrix. %, and propose a recursive algorithm to perform the construction.

\subsection{Non-Uniform Hypergraph}
A nonuniform hypergraph allows hyperedges to contain any number of vertices, relaxing the constraint of hyperedge uniformity. The adjacency structure for nonuniform hypergraphs can be represented as a series of adjacency tensors for each set of uniformly sized hyperedges.
% A \textit{undirected hypergraph} $\g=\{\Vs,\Es\}$ where $\Vs$ is a finite set and $\Es\subseteq \mathcal{P}(\Vs)\setminus \{\emptyset\}$, the power set of $\Vs$. The elements of $\Vs$ are called the nodes, and the elements of $\Es$ are called the hyperedges. A hypergraph is $k$-uniform if all hyperedges contain exactly $k$ vertices.

% \begin{defn}
% Let $\g= \{\Vs,\Es\}$ be a $k$-uniform hypergraph with $n=|\Vs|$ nodes. The \textit{adjacency tensor} $\tA\in\mathbb{R}^{n\times n\times \dots\times n}$ of $\g$ is a $k$-th order, $n$-dimensional, supersymmetric tensor is defined
% \begin{equation}\label{eq:98}
%     \tA_{j_1j_2\dots j_k} =
%     \begin{cases}
%         \frac{1}{(k-1)!} & \text{if $\{j_1,j_2,\dots,j_k\}\in \Es$}\\
%         0 &  \text{otherwise}\end{cases}.
% \end{equation}
% \end{defn}

\begin{defn}
Let $\g= \{\Vs,\Es\}$ be a non-uniform hypergraph with $n=|\Vs|$ nodes and maximum edge cardinality $K$. For each edge cardinality $k=2,\cdots,K$  we can associate an adjacency tensor
\begin{equation*}
\tA_k\in\mathbb{R}^{n\times n\times \dots\times n},
\end{equation*}
which is a $k$-th order, $n$-dimensional, supersymmetric tensor which captures all $k$-th order hyperedge interactions. Then $\g$ can be described by a collection of adjacency tensors $\{\tA_k\}_{k=2}^K$.
\end{defn}

\subsection{Nonuniform Hypergraph Dynamics with Outputs}

\begin{defn}
Given a $K$-nonuniform undirected hypergraph $\g$ with $n$ nodes, the dynamics of $\g$ with outputs $\yv\in \R^m$ is defined as
\begin{equation}\label{eq:hypdyn1non}
    \Sigma\begin{cases}
        \dot{\xv} & = \fv(\xv)=\sum_{k=2}^{K}\tA_{k}\xv^{k-1},\\
        \yv & = \gv(\xv)=\Cv\xv,
    \end{cases}
\end{equation}
where $\tA_k\in\mathbb{R}^{n\times n\times \dots \times n}$ is the adjacency tensors of $\g$, and $\Cv\in\mathbb{R}^{m\times n}$ is the output matrix.
\end{defn}

The nonuniform hypergraph dynamics (\ref{eq:hypdyn1non}) can be expressed equivalently as the unfolded tensor $\tA$ contracted with the Kronecker exponentiation of the state vector:
\begin{equation}\label{eq:hypdynkron1}
    \Sigma\begin{cases}
        \dot{\xv} & = \fv(\xv) = \sum_{k=2}^K \Av_k\xv^{[k-1]},\\
        \yv & = \gv(\xv) = \Cv\xv,
    \end{cases}
\end{equation}
% \begin{eqnarray}
% \dot{\xv}&=&\Av\xv^{[k-1]},\label{eq:hypdynkron1}\\
% \yv&=&\Cv\xv,\label{eq:hypdynkron2}
% \end{eqnarray}
where, $\Av_k=\tA_{k(p)}\in \R^{n\times n^{k-1}}$ is the $p$-th mode unfolding of $\tA_k$,{ and $\xv^{[k]}=\xv\otimes\xv\cdots\otimes\xv$ with $\otimes$ being the Kronecker product}. Since $\tA_k$ is a super-symmetric tensor, all $p$-th mode unfoldings give rise to the same matrix.

 Tensor unfoldings:
Tensor unfolding is a fundamental operation in tensor computations \cite{kolda2009tensor}. In order to unfold a tensor $\tT\in\R^{n_1\times\cdots \times n_k}$ into a vector or a matrix, we use an index mapping function $ivec(\cdot,\mathcal{N}):\Z\times\Z\times \stackrel{\scriptscriptstyle k}{\cdots} \times \Z \rightarrow \Z$, which is defined as
 \begin{equation*}
 ivec(\mathcal{J},\mathcal{N}) = j_1+\sum_{i=2}^k(j_i-1)\prod_{l=1}^{i-1}n_l,
\end{equation*}
 where $\mathcal{J}=\{j_1,j_2,\cdots,j_k\}$ and $\mathcal{N}=\{n_1,\cdots,n_k\}$ are the sets of indices and mode sizes of $\tT$, respectively. The function $ivec$ is an injective map from the set of the tensor indices $\{j_1,j_2,\cdots,j_k\}$ to a unique vector index and is used to unfold or matricize tensors.

 Tensor $p$-mode unfolding \cite{kolda2009tensor}:
    Given a $k$th-order tensor $\tT\in \mathbb{R}^{n_1\times  \dots \times n_k}$, the $p$-mode unfolding of $\tT$, denoted by $\tT_{(p)}\in\R^{n_p \times (n_1\cdots n_{p-1}n_{p+1}\cdots n_k)}$, is defined as
 \begin{equation*}
 (\tT_{(p)})_{ij}=\tT_{j_1\cdots j_{p-1}ij_{p+1}\cdots j_k},
\end{equation*}
 where $j=ivec(\tilde{\mathcal{J}}, \tilde{\mathcal{N})}$ with $\tilde{\mathcal{J}}=\mathcal{J}\backslash \{j_p\}$ and $\tilde{\mathcal{N}}=\mathcal{N}\backslash \{n_p\}$.

\subsection{Observability Criterion}
To determine the nonlinear observabilty matrix (\ref{eq:obsrmat}) for the systems (\ref{eq:hypdyn1non}) and (\ref{eq:hypdynkron1}), we compute the Lie derivatives of the system output along the flow $\fv$ of the system state, which can be decomposed as:
\begin{equation*}
\fv(\xv) = \sum_{k=2}^K \fv_k(\xv),
\end{equation*}
where, $\fv_k=\Av_k\xv^{[k-1]}$. Also note that
\begin{equation*}
 L_{\fv}\gv(\xv)=L_{\sum_{k_1=2}^K \fv_{k_1}}\gv(\xv)=\sum_{k_1=2}^{K} L_{\fv_{k_1}}\gv(\xv).
\end{equation*}
and
\begin{equation*}
 L^2_{\fv}\gv(\xv)=\sum_{k_1=2}^{K} \sum_{k_2=2}^K  L_{\fv_{k_2}}(L_{\fv_{k_1}}\gv(\xv)),
\end{equation*}
and, hence,
\begin{equation*}
 L^n_{\fv}\gv(\xv)=\sum_{k_1=2}^{K} \sum_{k_2=2}^K \cdots \sum_{k_n=2}^K L_{\fv_{k_n}}\cdots L_{\fv_{k_2}}L_{\fv_{k_1}}\gv(\xv).
\end{equation*}

For each $\fv_k$ one can show (time derivatives here are along the appropriate vector fields),
\begin{eqnarray*}
    L_{\fv_k}^0\gv(\xv) & =& \Cv\xv,\\
    L_{\fv_k}^1\gv(\xv) & = & \frac{d}{dt}\Cv\xv=\Cv\Av_k\xv^{[k-1]},\\
    L_{\fv_k}^2\gv(\xv) & = & \frac{d}{dt}\Cv\Av_k\xv^{[k-1]}\\
              & = & \Cv\Av_k\overline{\Bv}_{k2} \xv^{2k-3},\\
    & \vdots &\\
    L_{\fv_k}^n\gv(\xv) & = & \Cv\Av_k\overline{\Bv}_{k2}\dots \overline{\Bv}_{kn} \xv^{[nk-(2n-1)]} \quad \forall n>2,
\end{eqnarray*}
where $\overline{\Bv}_{kp}$ is given by
\begin{equation}\label{eq:B}
\overline{\Bv}_{kp}=\sum_{i=1}^{(p-1)k-(2p-3)} \overbrace{\Id\otimes\dots\otimes \overbrace{\Av_k}^{i\text{-th pos.}}\otimes\dots\otimes \Id}^{(p-1)k-(2p-3) \text{times}}.
\end{equation}
On the other hand, for mixed terms,
\begin{eqnarray*}
L_{\fv_{k_2}}L_{\fv_{k_1}}\gv(\xv) & = & \frac{d}{dt}\Cv\Av_{k_1}\xv^{[{k_1}-1]}\notag \\
& = & \Cv\Av_{k_1}\frac{d}{dt}\bigg(\overbrace{\xv\otimes \xv\dots\otimes\xv}^{\text{$k_1-1$ times}}\bigg)\\
              & = & \Cv\Av_{k_1}\bigg(\dot{\xv}\otimes\dots\otimes \xv+\cdots+\xv\otimes\dots\otimes\dot{\xv}\bigg)\\ %\xv\otimes\dot{\xv}\otimes\dots\otimes \xv
              & = & \Cv\Av_{k_1}\bigg(\sum \xv\otimes\dots\otimes \Av_{k_2}\xv^{[k_2-1]}\otimes\dots\otimes \xv\bigg)\\
              & = & \Cv\Av_{k_1}\bigg[\bigg(\sum \Id\otimes\dots \Av_{k_2}\otimes\dots \Id\bigg) \xv^{[k_2+k_1-3]}\bigg]\nonumber\\
              & = & \Cv\Av_{k_1}\Bv_{k_2k_1} \xv^{[k_1+k_2-3]},\\
\end{eqnarray*}
where,
\begin{equation}\label{eq:Bk1}
\Bv_{k_2k_1}=\sum_{i=1}^{k_1-1} \overbrace{\Id\otimes\dots\otimes \overbrace{\Av_{k_2}}^{i\text{-th pos.}}\otimes\dots\otimes \Id}^{k_1-1 \text{times}}.
\end{equation}
Similarly,
\begin{eqnarray*}
L_{\fv_{k_3}}L_{\fv_{k_2}}L_{\fv_{k_1}}\gv(\xv)  &=& \Cv\Av_{k_1}\Bv_{k_2k_1}\Bv_{k_3k_2k_1}  \xv^{[k_1+k_2+k_3-5]},
\end{eqnarray*}
where,
\begin{equation}\label{eq:Bk1k2}
\Bv_{k_3k_2k_1}=\sum_{i=1}^{k_1+k_2-3} \overbrace{\Id\otimes\dots\otimes \overbrace{\Av_{k_3}}^{i\text{-th pos.}}\otimes\dots\otimes \Id}^{k_1+k_2-3 \text{times}}.
\end{equation}
Thus, in general,
\begin{eqnarray*}
&&L_{\fv_{k_n}}\cdots L_{\fv_{k_3}}L_{\fv_{k_2}}L_{\fv_{k_1}}\gv(\xv)\\
&=& \Cv\Av_{k_1}\Bv_{k_2k_1}\Bv_{k_3k_2k_1} \cdots \Bv_{k_n\cdots k_3k_2k_1}  \xv^{[\sum_{i=1}^n k_{i}-(2n-1)]},
\end{eqnarray*}
with,
\begin{equation*}\label{eq:Bk1k2kn}
\Bv_{k_n\cdots k_3k_2k_1}=\sum_{i=1}^{\sum_{i=1}^{n-1}k_i-(2n-3)} \overbrace{\Id\otimes\dots\otimes \overbrace{\Av_{k_n}}^{i\text{-th pos.}}\otimes\dots\otimes \Id}^{\sum_{i=1}^{n-1}k_i-(2n-3) \text{times}}.
\end{equation*}
Note that if $n=p$ and $k_1=k_2=\cdots k_p=k$, i.e., the uniform hypergraph case, then
\begin{eqnarray*}\label{eq:Bkp}
\Bv_{\underbrace{k\cdots kkk}_{p-\text {times}}}&=&\sum_{i=1}^{\sum_{i=1}^{p-1}k-(2p-3)} \overbrace{\Id\otimes\dots\otimes \overbrace{\Av_{k}}^{i\text{-th pos.}}\otimes\dots\otimes \Id}^{\sum_{i=1}^{p-1}k-(2p-3) \text{times}}\\
&=&\sum_{i=1}^{(p-1)k-(2p-3)} \overbrace{\Id\otimes\dots\otimes \overbrace{\Av_{k}}^{i\text{-th pos.}}\otimes\dots\otimes \Id}^{(p-1)k-(2p-3) \text{times}},\nonumber
\end{eqnarray*}
% \textcolor{blue}
{reduces to the matrix (see \cite{PiSuBlRa2023}) which arise in calculating NOM for uniform hypergraph case as expected}.

Finally,
\begin{align*}
 &L^n_{\fv}\gv(\xv)=\sum_{k_1=2}^{K} \sum_{k_2=2}^K \cdots \sum_{k_n=2}^K \Cv\Av_{k_1}\\
 &\Bv_{k_2k_1}\Bv_{k_3k_2k_1} \cdots \Bv_{k_n\cdots k_3k_2k_1} \xv^{[\sum_{i=1}^n k_{i}-(2n-1)]} ,
\end{align*}
and we can express the NOM as,
\begin{equation}\label{eq:hgobsv}
    \mathcal{O}(\xv)=\nabla_{\xv}\begin{pmatrix}
    \Cv\xv\\
    \sum_{k_1=2}^K\Cv\Av_{k_1}\xv^{k_1-1}\\
    \sum_{k_1=2}^{K} \sum_{k_2=2}^K \Cv\Av_{k_1}\Bv_{k_2k_1}\xv^{[k_1+k_2-3]}  \\
    \\
    \vdots\\
    L^n_{\fv}\gv
    %\sum_{k_1=2}^{K} \sum_{k_2=2}^K \cdots \sum_{k_n=2}^K \Cv\Av_{k_1}\Bv_{k_2k_1}\Bv_{k_3k_2k_1} \cdots \Bv_{k_n\cdots k_3k_2k_1} \xv^{[\sum_{i=1}^n k_{i}-(2n-1)]}
    \end{pmatrix},
\end{equation}
where, we have used $r=n$ as per the Remark \ref{remrn}. From Theorem \ref{thm1}, when $\mbox{rank}(\mathcal{O}(\xv))=\dim(\xv),$ systems (\ref{eq:hypdyn1non}) and equivalently (\ref{eq:hypdynkron1}) are observable.

\begin{figure*}[t!]
    \centering
    \includegraphics[width=\textwidth]{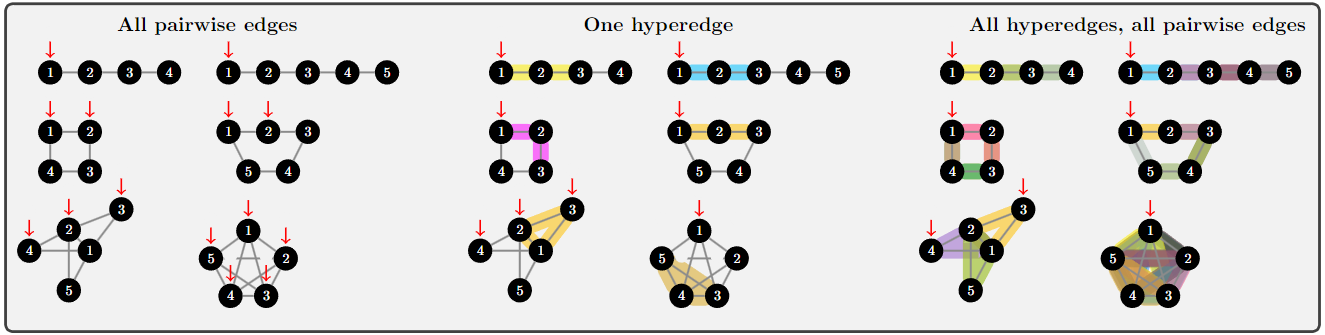}
    \caption{\textbf{(Left)} Hypergraphs containing only pairwise edges are shown. Red arrows indicate a MON (minimum observable note) set to render each structure localy weakly observable. The top row contains two hyperchains, the middle row is two hyperrings, the bottom row contains the clique expansion of a uniform hyperstar and a complete hypergraph. \textbf{(Middle)} A single 3-way hyperedge has been added to each of the 2-uniform hypergraphs on the left, improving the observability of all hypergraphs. \textbf{(Right)} All 3-way hyperedges have been added to the hyperchain, hyperring, hyperstar and complete hypergraphs while still including the pairwise structure. c.f \cite{PiSuBlRa2023} for the MON on similar 3-uniform hypergraphs.}
    \label{fig:enter-label}
\end{figure*}

\color{black}
\begin{prop}
    Given a nonuniform hypergraph $\g$ with a maximum hyperedge cardinality of $k,$ if $\g$ is locally weakly observable when all hyperedges with cardinality less than $k$ are removed, then $\g$ is locally weakly observable.
\end{prop}

\textit{Proof.} Consider a $k$-uniform hypergraph $\g$ that is locally weakly observable such that the symbolic NOM $\mathcal{O}(\xv)$ is of full rank. % Suppose a hyperedge of cardinality less than $k$ were added to $\h$
The inclusion of a hyperedge with cardinality less than $k$
in $\g$ will introduce polynomial terms of $\xv$ into the symbolic NOM; however, the highest order polynomial terms in the NOM are generated by hyperedges of maximum cardinality, such that lower order polynomial terms cannot decrease the symbolic rank. Since hyperedges of lower cardinality introduce lower order polynomial terms into $\mathcal{O}(\xv),$ the inclusion of hyperedges with fewer than $k$ vertices cannot decrease the symbolic rank of the hypergraph with only $k$-way hyperedges.
\hfill $\blacksquare$

% \textit{Proof.}
%     If the uniform hypergraph with the highest cardinality hyperedges from $\h$ is observable, then the symbolic matrix $\Ov(\xv)$ is full rank. Moreover, the highest order polynomial terms of $\Ov(\xv)$ are generated with the highest cardinality hyperedges. All hyperedges with lower cardinality will generate lower order polynomial terms, which will never decrease the rank of $\mathcal{O}(\xv).$ \hfill \blacksquare

In Fig. \ref{fig:enter-label}, several hypergraph topologies %(chains, rings, stars, and complete graphs)
are drawn with the minimum observable nodes (MON) required to make the hypergraphs observable. The number of 2-uniform hyperedges is fixed, and the number of 3-uniform hyperedges varies. %As the number of hyperedges increases, several MON sets may be chosen, and 
We utilized the greedy selection procedure of \cite{PiSuBlRa2023}, where only 3-uniform hypergraphs were considered, to draw Fig \ref{fig:enter-label}. Consistent with Proposition 15, the increase in hyperedges uniformly decreases the size of the MON, which further motivates the investigation of higher order observability.

% In Figure \ref {fig:enter-label} we depict various hypergraphs and in particular the minimum number of observable nodes (MON) needed to achieve of observability of such systems.  These can be calculated algorithmically using symbolic computation as in \cite{PiSuBlRa2023}. \textcolor{blue}{I think this section requires bit details, may be Joshua can you address that?}
\color{black}

\section{Conclusion}
We have analyzed observability for non-uniform hypergraphs, applying the theory developed earlier for the uniform 
case, and illustrated the work on several low-dimensional examples. Future work will include applying these 
ideas to large data sets, particularly those arising from biological applications. 

{\bf Acknowledgement:} We would like to thank the reviewers for very helpful comments regarding the exposition. 

\bibliographystyle{johd}
\bibliography{bib}

\end{document}